\documentclass[oneside,english]{amsart}
\usepackage{mathptmx}
\usepackage[T1]{fontenc}
\usepackage[latin9]{inputenc}
\usepackage{amsthm}
\usepackage{graphicx}
\usepackage{amssymb}

\newcommand{\lyxline}[1][1pt]{%
  \par\noindent%
  \rule[.5ex]{\linewidth}{#1}\par}

\numberwithin{equation}{section} 
\numberwithin{figure}{section} 
\theoremstyle{plain}
\theoremstyle{plain}
\newtheorem{thm}{Theorem}
  \theoremstyle{plain}
  \newtheorem{cor}[thm]{Corollary}

\usepackage{babel}

\begin{document}

\title{On Plouffe's Ramanujan Identities}

\author{Linas Vepštas }

\email{<linasvepstas@gmail.com>}

\date{14 May 2006 Revised 27 November 2010}
\begin{abstract}
Recently, Simon Plouffe has discovered a number of identities for
the Riemann zeta function at odd integer values. These identities
are obtained numerically and are inspired by a prototypical series
for Apéry's constant given by Ramanujan:\[
\zeta(3)=\frac{7\pi^{3}}{180}-2\sum_{n=1}^{\infty}\frac{1}{n^{3}\left(e^{2\pi n}-1\right)}\]
Such sums follow from a general relation given by Ramanujan, which
is rediscovered and proved here using complex analytic techniques.
The general relation is used to derive many of Plouffe's identities
as corollaries. The resemblance of the general relation to the structure
of theta functions and modular forms is briefly sketched.
\end{abstract}

\subjclass[2000]{11M06}

\maketitle

\section{Introduction}

Inspired by an identity for $\zeta(3)$ given in Ramanujan's notebooks
\cite[chapter 14, formulas 25.1 and 25.3]{Ber-II},\[
\zeta(3)=\frac{7\pi^{3}}{180}-2\sum_{n=1}^{\infty}\frac{1}{n^{3}\left(e^{2\pi n}-1\right)}\]
Plouffe describes a set of similar identities\cite{Plo98},\cite{Plo06}
that were discovered numerically using arbitrary-precision software.
For example, Plouffe gives an identity for $\zeta(7)$: \[
\zeta(7)=\frac{19\pi^{7}}{56700}-2\sum_{n=1}^{\infty}\frac{1}{n^{7}\left(e^{2\pi n}-1\right)}\]
This text provides an analytically derived formula for expressions
of this type. The resulting general formula, valid for integer $m\ge1$,
is \begin{eqnarray*}
\zeta(4m-1) & = & -2\sum_{n=1}^{\infty}\frac{1}{n^{4m-1}\left(e^{2\pi n}-1\right)}\\
 & - & \frac{1}{2}\;(2\pi)^{4m-1}\;\sum_{j=0}^{2m}\left(-1\right)^{j}\frac{B_{2j}}{(2j)!}\frac{B_{4m-2j}}{(4m-2j)!}\end{eqnarray*}
where $B_{k}$ is the $k$'th Bernoulli number. The above is a special
case of a yet more general formula, derived and presented in a later
section, allowing pairs of such sums to be related. From this, one
may obtain expressions such as\[
\zeta(3)=\frac{37\pi^{3}}{900}-\frac{2}{5}\sum_{n=1}^{\infty}\frac{1}{n^{3}}\left[\frac{4}{e^{\pi n}-1}+\frac{1}{e^{4\pi n}-1}\right]\]
There are an (uncountable) infinity of similar sums, each giving a
different series summation for $\zeta(4m-1)$. Taking linear combinations
of these, one may choose to cancel the zeta terms, to obtain summations
for odd powers of $\pi$. Thus, for example, combining the above with
Ramanujan's series for Apéry's constant, one gets:\[
\frac{\pi^{3}}{180}=\sum_{n=1}^{\infty}\frac{1}{n^{3}}\left[\frac{4}{e^{\pi n}-1}-\frac{5}{e^{2\pi n}-1}+\frac{1}{e^{4\pi n}-1}\right]\]
Again, there are an uncountable infinity of such relations. 

Plouffe also notes similar relations for the other odd integers; for
example, \[
\zeta(5)=\frac{\pi^{5}}{294}-\frac{72}{35}\sum_{n=1}^{\infty}\frac{1}{n^{5}\left(e^{2\pi n}-1\right)}-\frac{2}{35}\sum_{n=1}^{\infty}\frac{1}{n^{5}\left(e^{2\pi n}+1\right)}\]
The general form for this type of expression may be shown to be \begin{gather*}
\left[1+(-4)^{m}-2^{4m+1}\right]\zeta(4m+1)\qquad\qquad\qquad\qquad\qquad\\
\quad=2\sum_{n=1}^{\infty}\frac{1}{n^{4m+1}\left(e^{2\pi n}+1\right)}\qquad\qquad\qquad\qquad\qquad\\
\qquad+2\left[2^{4m+1}-(-4)^{m}\right]\sum_{n=1}^{\infty}\frac{1}{n^{4m+1}\left(e^{2\pi n}-1\right)}\qquad\\
\qquad\qquad+(2\pi)^{4m+1}\sum_{j=0}^{m}(-4)^{m+j}\frac{B_{4m-4j+2}}{(4m-4j+2)!}\frac{B_{4j}}{(4j)!}\\
\qquad\qquad\qquad+\frac{1}{2}(2\pi)^{4m+1}\;\sum_{j=0}^{2m+1}(-4)^{j}\frac{B_{4m-2j+2}}{(4m-2j+2)!}\frac{B_{2j}}{(2j)!}\end{gather*}

The methods described in this text also allow for a large generalization
of these types of sums. Defining \[
P_{k}(\tau)=\sum_{n=1}^{\infty}\frac{1}{n^{k}\left(e^{2\pi in\tau}-1\right)}\]
these generalizations follow from a modular equation relating $P_{k}(\tau)$
to $P_{k}(-1/\tau)$ for odd integers $k$, the derivation and proof
of which is the one of the main topics of this note. The modular relation
is not new; it appears in Ramanujan's Notebooks \cite[Chapter 14 Entry 21]{Ber-II}
as \begin{multline*}
\alpha^{-n}\left\{ \frac{1}{2}\zeta\left(2n+1\right)+\sum_{k=1}^{\infty}\frac{k^{2n-1}}{e^{2\alpha k}-1}\right\} \\
\qquad=(-\beta)^{-n}\left\{ \frac{1}{2}\zeta\left(2n+1\right)+\sum_{k=1}^{\infty}\frac{k^{2n-1}}{e^{2\beta k}-1}\right\} \\
\qquad\qquad-2^{2n}\sum_{k=0}^{n+1}(-1)^{k}\frac{B_{2k}}{(2k)!}\frac{B_{2n+2-2k}}{(2n+2-2k)!}\alpha^{n+1-k}\beta^{k}\end{multline*}
where $\alpha>0$, $\beta>0$ with $\alpha\beta=\pi^{2}$ and $n$
any positive integer. Berndt implies that this formula is the most
studied of all the notebooks; it has been independently discovered
perhaps a half-dozen times, and proven twice as often. It has been
generalized to $L$-functions, and to rational values of $k$; Berndt
provides a long list\cite{Ber-II} of the various proofs and generalizations
made. 

Much of this paper is devoted to (yet another! independently discovered)
proof of this relation, followed by a series of lemmas that provide
the connection to Plouffe's results. In searching for curious and
interesting special cases of this relation, one senses that only the
tip of the iceberg has been seen. Unexplored possibilities include,
for example, considering $\tau\in\mathbb{Q}[i]$, the field of Gaussian
rationals, or from considering Diophantine roots of quadratics. 

The rest of this paper is roughly laid out as follows: The second
section provides a review of previous related results. The third section
gives a relationship between the sums and the polylogarithm, and thence
to an integral on the complex plane. The fourth section examines the
related contour integral, which is easily integrated via Cauchy's
residue theorem to give a finite sum involving the Bernoulli numbers.
The fifth section relates the contour integral to the polylogarithm
integral, thus resulting in a functional equation for $P_{k}(\tau)$.
The sixth section applies the functional equation, providing a variety
of lemmas, many of which explain Plouffe's discoveries. The seventh
section give a pair of relationships on the Bernoulli numbers that
arise naturally in this context. The eighth section explores the modular
nature of the relations on $P_{k}(\tau)$, followed by a conclusion.
An appendix gives a derivation of an integral representation of the
polylogarithm, that is central to the analysis.

\section{Related sums}

A large number of similar sums have been explored before; this section
reviews some of these. Perhaps the most forthright is a sum given
by Ramanujan in a famous letter to Hardy\cite{Har21}, stating that
\[
\frac{1^{13}}{e^{2\pi}-1}+\frac{2^{13}}{e^{4\pi}-1}+\frac{3^{13}}{e^{6\pi}-1}+\cdots=\frac{1}{24}\]
A generalization of this sum, \[
\sum_{n=1}^{\infty}\frac{n^{4k+1}}{e^{2\pi n}-1}=\frac{B_{4k+2}}{4(2k+1)}\]
is proved by Berndt\cite{Ber77}, and attributed to Glaisher\cite{Gla1889}.
This and many related results are derived by Zucker\cite{Zuc79},
based on the theory of Jacobian elliptic functions. A similar result
is stated by Apostol in the form of an exercise\cite[see exercise 15 at end of chapter 1.]{Apo90}:\[
\sum_{n=1;n\mbox{ odd}}^{\infty}\frac{n^{4k+1}}{1+e^{n\pi}}=\frac{2^{4k+1}-1}{8k+4}B_{4k+2}\]

Many sums resembling those in this note are given by Zucker\cite{Zuc84}.
Some of these are \[
\sum_{n=1}^{\infty}\frac{1}{n\left(e^{2\pi nx}-1\right)}=\frac{1}{2}\sum_{m=1}^{\infty}\frac{\coth(\pi mx)-1}{m}\]
 \[
\sum_{n=1}^{\infty}\frac{(-1)^{n}}{n\left(e^{2\pi nx}-1\right)}=\frac{1}{2}\sum_{m=1}^{\infty}(-1)^{m}\,\frac{\coth(\pi mx)-1}{m}\]
 \[
\sum_{n=1}^{\infty}\frac{1}{n\left(e^{2\pi nx}+1\right)}=\frac{1}{2}\sum_{m=1}^{\infty}\frac{1-\tanh(\pi mx)}{m}\]
\[
\sum_{n=1}^{\infty}\frac{(-1)^{n}}{n\left(e^{2\pi nx}+1\right)}=\frac{1}{2}\sum_{m=1}^{\infty}(-1)^{m}\;\frac{1-\tanh(\pi mx)}{m}\]
The fount of inspiration for such sums is Ramanujan. Sums given in
Chapter 14 of Part II of the Ramanujan's Notebooks\cite{Ber-II} include
entry 8:\[
\alpha\sum_{n=1}^{\infty}\frac{\sinh(2\alpha nk)}{e^{2\alpha^{2}n}-1}+\beta\sum_{n=1}^{\infty}\frac{\sin(2\beta nk)}{e^{2\beta^{2}n}-1}=\frac{\alpha}{4}\coth(\alpha k)-\frac{\beta}{4}\cot(\beta k)-\frac{k}{2}\]
and another similar one relating cos and cosh. Above, $k$ is any
positive integer, $\alpha\beta=\pi$ and $0<\beta k<\pi$. 

Entry 13 generalizes the sums mentioned previously, \[
\alpha^{k}\sum_{n=1}^{\infty}\frac{n^{2k-1}}{e^{2\alpha n}-1}-(-\beta)^{k}\sum_{n=1}^{\infty}\frac{n^{2k-1}}{e^{2\beta n}-1}=\left[\alpha^{k}-(-\beta)^{k}\right]\frac{B_{2k}}{4k}\]
This time, one takes $\alpha\beta=\pi^{2}$and $k>1$ an integer.
The above follows from sums on the divisor function, as is frequently
noted. 

Sums involving pairs of Bernoulli numbers also appear in the analysis
of the Dedekind eta function. Thus, Sho Iseki's transformation formula,
as described by Apostol\cite[see theorem 3.5]{Apo90}, is \[
\Lambda(\alpha,\beta,z)=\Lambda\left(1-\beta,\alpha,\frac{1}{z}\right)-\pi z\sum_{n=0}^{2}\left(\begin{array}{c}
2\\
n\end{array}\right)\frac{B_{2-n}(\alpha)B_{n}(\beta)}{(iz)^{n}}\]
where $\Lambda$ is given by \[
\Lambda(\alpha,\beta,z)=\sum_{r=0}^{\infty}\left[\lambda(z(r+\alpha)-i\beta)+\lambda(z(r+1-\alpha)+i\beta)\right]\]
and \[
\lambda(x)=\sum_{m=1}^{\infty}\frac{e^{-2\pi mx}}{m}\]
A sum linking the Bernoulli and Euler numbers is given by Berndt\cite{Ber78}:
\begin{multline*}
\sum_{n=1}^{\infty}(-1)^{n+1}\frac{\mbox{sech}\left[(2n-1)\pi\sqrt{3}/2\right]}{(2n-1)^{6k+1}}\\
\qquad=\frac{1}{2}(-1)^{k+1}\pi^{6k+1}\sum_{m=0}^{3k}\frac{E_{2m+1}}{(2m+1)!}\,\frac{B_{6k-2m}}{(6k-2m)!}\cos\left[(2m+1)\frac{\pi}{3}\right]\end{multline*}

Perhaps the results that are closest to those presented in this paper
are those noted by Borwein, \emph{et al} \cite[Section 5]{Bor00},
and in particular, giving a very similar result involving $\zeta(4m-1)$
and $\zeta(4m+1)$.

\section{The Polylogarithm}

The recurring theme in Plouffe's identities is the sum \[
S_{s}(x)=\sum_{n=1}^{\infty}\frac{1}{n^{s}\left(e^{xn}-1\right)}\]
with $s$ usually an odd positive integer and $x=\pi$ or $x=2\pi$
or possibly other interesting values, such as $x=\pi\sqrt{m}$ for
some integer $m$. This sum may be converted into a sum over polylogarithms,
and subsequently into an integral. The integral may, after some difficulties,
be converted into a contour integral, whereupon it may be evaluated
by Cauchy's residue theorem. The result is a finite sum whose general
structure resembles those of Plouffe's and Ramanujan's identities.
This section develops the first part of this analysis. 

To find the polylogarithm, one expands\begin{eqnarray*}
S_{s}(x) & = & \sum_{n=1}^{\infty}\sum_{m=0}^{\infty}\frac{e^{-xn(m+1)}}{n^{s}}\\
 & = & \sum_{m=1}^{\infty}{\rm Li}_{s}\left(e^{-xm}\right)\end{eqnarray*}
which is generally valid for $\Re x>0$. The above is easily obtained
by applying the expansion \[
\frac{1}{1-z}=\sum_{m=0}^{\infty}z^{m}\]
and the series definition of the polylogarithm:\[
{\rm Li}_{s}(z)=\sum_{n=1}^{\infty}\frac{z^{n}}{n^{s}}\]
The polylogarithm may be expressed in terms of an integral as\[
{\rm Li}_{s}\left(e^{-u}\right)=\frac{1}{2\pi i}\int_{c-i\infty}^{c+i\infty}\Gamma(z)\zeta(z+s)u^{-z}\, dz\]
A derivation of this is given in the appendix. Here, $\Gamma(z)=(z-1)!$
is the classical Gamma function. The line of integration is taken
to be to the right of all of the poles in the integrand, namely, $c>1$.
Using this in the summation, one obtains\begin{eqnarray*}
S_{s}(x) & = & \frac{1}{2\pi i}\sum_{m=1}^{\infty}\int_{c-i\infty}^{c+i\infty}\Gamma(z)\zeta(z+s)\left(xm\right)^{-z}\, dz\\
 & = & \frac{1}{2\pi i}\int_{c-i\infty}^{c+i\infty}\frac{\Gamma(z)}{x^{z}}\zeta(z+s)\zeta(z)\, dz\end{eqnarray*}
The exchange of the order of summation and integration is justified
precisely when one has $c>1$. The last integral has poles at $z=1$
and $z+s=1$ coming from the zeta functions and poles at all of the
non-positive integers coming from the Gamma function. The last integral
shows that the series is the inverse Mellin transform of $\Gamma(z)\zeta(z+s)\zeta(z)$.

If the integral can somehow be converted into a closed contour on
the left, then it may be evaluated in a straight-forward way by means
of Cauchy's residue theorem. Performing this closure is in fact harder
than one might hope, as there are non-zero contributions to the contour
from its closure. The next section evaluates the Cauchy integral,
assuming that the contour can be closed. The section after that computes
the contributions from closing the contour integral. Upon doing this,
Plouffe's identities, and many more, become available.

\section{The Contour Integral}

Define the contour integral as \[
I_{s}(x)=\frac{1}{2\pi i}\oint_{\gamma}\frac{\Gamma(z)}{x^{z}}\zeta(z+s)\zeta(z)dz\]
where the contour $\gamma$ encircles the poles at $z=1$, $z+s=1$
and $z=0,-1,-2,\ldots$ in the usual, right-handed fashion. Then,
one uses Cauchy's theorem, which states that \[
f(a)=\frac{1}{2\pi i}\oint\,\frac{f(z)}{z-a}\, dz\]
for simple poles, and that\[
f^{\,\prime}(a)=\frac{1}{2\pi i}\oint\,\frac{f(z)}{\left(z-a\right)^{2}}\, dz\]
for double poles. For the pole at $z=1$, one obtains the residue\[
{\rm Res}(z=1)=\frac{\zeta(s+1)}{x}\]
For the poles at $z=-n$, one obtains the residue\[
{\rm Res}(z=-n)=\frac{(-x)^{n}}{n!}\zeta(s-n)\zeta(-n)\]
and so one has\[
I_{s}(x)={\rm Res}(z=1-s)+\frac{\zeta(s+1)}{x}+\sum_{n=0}^{\infty}\frac{(-x)^{n}}{n!}\zeta(s-n)\zeta(-n)\]
For $s$ not an integer, one has \[
{\rm Res}(z=1-s)=\frac{\Gamma(1-s)}{x^{1-s}}\zeta(1-s)\]
However, the interesting case is for $s=k$ a positive integer. In
this case, the pole overlays another pole from the Gamma, and one
has a double pole. This is just a little trickier to evaluate:\begin{eqnarray*}
\frac{1}{2\pi i}\oint_{z=1-k}\frac{\Gamma(z)}{x^{z}}\zeta(z+s)\zeta(z)dz & = & \frac{1}{2\pi i}\oint_{z=1-k}\frac{f(z)}{(z+k-1)^{2}}dz\\
 & = & \left.\frac{d}{dz}\left[(z+k-1)^{2}\Gamma(z)\;\zeta(z+k)\frac{\zeta(z)}{x^{z}}\right]\right|_{z=1-k}\end{eqnarray*}
To perform the derivative, one will need to use the identities\[
\left.\frac{d}{ds}(s-1)\zeta(s)\right|_{s=1}=\gamma\]
where $\gamma=0.577\ldots$ is the Euler-Mascheroni constant, and
\[
\left.\frac{d}{dz}(z+n)\Gamma(z)\right|_{z=-n}=(-1)^{n}\,\frac{\psi(n+1)}{\Gamma(n+1)}=(-1)^{n}\,\frac{H_{n}-\gamma}{n!}\]
where $\psi(z)$ is the digamma function, and $H_{n}$ is the $n$'th
harmonic number. Putting these together, one obtains \[
\frac{1}{2\pi i}\oint_{z=1-k}\frac{\Gamma(z)}{x^{z}}\zeta(z+s)\zeta(z)dz=\frac{(-x)^{k-1}}{(k-1)!}\left[\zeta^{\prime}(1-k)+\left(H_{k-1}-\ln2\pi\right)\zeta(1-k)\right]\]
Adding this to the other contributions, one gets

\begin{eqnarray*}
I_{k}(x) & = & \frac{\zeta(k+1)}{x}+\!\sum_{\begin{array}{c}
n=0\\
n\ne k-1\end{array}}^{\infty}\!\frac{(-x)^{n}}{n!}\zeta(k-n)\zeta(-n)\\
 &  & +\,\frac{(-x)^{k-1}}{(k-1)!}\left[\zeta^{\prime}(1-k)+\left(H_{k-1}-\ln2\pi\right)\zeta(1-k)\right]\end{eqnarray*}
When $k$ is an odd integer, the above simplifies in two ways. First,
$\zeta(1-k)$ vanishes, because the zeta function vanishes at all
negative even integers. Similarly, the infinite sum becomes finite:
when $k$ is an odd integer, one has either $\zeta(k-n)=0$ or $\zeta(-n)=0$
for all $n>k$. Thus, for $k$ an odd integer, one has \[
I_{k}(x)=\frac{\zeta(k+1)}{x}+\!\sum_{\begin{array}{c}
n=0\\
n\ne k-1\end{array}}^{k}\!\frac{(-x)^{n}}{n!}\zeta(k-n)\zeta(-n)+\frac{(-x)^{k-1}}{(k-1)!}\zeta^{\prime}(1-k)\]
For the remainder of the paper, it is assumed that, in this context,
$k>1$ is an odd integer, unless explicitly stated otherwise. The
above evaluation of the contour integral is the main result of this
section. To see that this is a key result, one may substitute $k=3$
and $x=2\pi$ to obtain \[
I_{3}(2\pi)=-\zeta(3)+\frac{7\pi^{3}}{180}\]
which should be recognizable as a portion of Ramanujan's identity.
For $k=7$, one has \[
I_{7}(2\pi)=-\zeta(7)+\frac{19\pi^{7}}{56700}\]
which resembles one of the results given by Plouffe. To complete the
connection, one must relate the contour integral $I_{k}(x)$ to the
sum $S_{k}(x)$. This is done in the next section. 

First, however, to drive the point home, one must observe that most
of the terms in the above expression are rational multiples of powers
of $\pi$. This follows from the zeta function being related to the
Bernoulli numbers $B_{n}$ at even integers:\[
\zeta(2n)=(-1)^{n+1}\frac{(2\pi)^{2n}B_{2n}}{2(2n)!}\]
for integer $n\ge0$. At the negative values, one has\[
\zeta(-n)=-\frac{B_{n+1}}{n+1}\]
while the derivative is \[
\zeta^{\prime}(-2n)=(-1)^{n}\frac{(2n)!}{2(2\pi)^{2n}}\zeta(2n+1)\]
for integer $n>1$. Using these in the above expression for $I_{k}(x)$,
and rearranging terms a bit, one gets \begin{eqnarray*}
-2I_{k}(x) & = & \zeta(k)\left[1-\left(\frac{x}{2\pi i}\right)^{k-1}\right]\\
 & + & \frac{1}{x\;(k+1)!}\;\sum_{j=0}^{(k+1)/2}\left(\!\begin{array}{c}
k+1\\
2j\end{array}\!\right)x^{2j}\left(2\pi i\right)^{k+1-2j}B_{2j}B_{k+1-2j}\end{eqnarray*}
The above introduces the binomial coefficient\[
\left(\!\begin{array}{c}
n\\
k\end{array}\!\right)=\frac{n!}{k!(n-k)!}\]
Although the imaginary number $i=\sqrt{-1}$ appears in the above,
it is always squared, and thus is just a sign-keeping device. Every
term in the sum is purely real. 

When the contour integral is written in this form, it may now be seen
that for $x$ being any rational multiple of $\pi$, that is, $x=p\pi/q$
for any integers $p,q$, that the coefficient of $\zeta(k)$ is a
rational number, and that the second term is another rational times
$\pi^{k}$. 

A further curiosity in this regard is noted by Plouffe: if one takes
$x=\pi\sqrt{p/q}$ for integers $p$ and $q$, one also gets simple
expressions: because $k$ is odd, the coefficient of $\zeta(k)$ is
still a rational, and the coefficient of $\pi^{k}$ is $\sqrt{p/q}$
times some rational. For rational $x$, one still has that the sum
is a rational polynomial in $\pi^{2}$, and for $x=\sqrt{\pi p/q}$,
one still has that the sum is a rational polynomial in $\pi$. The
ocean-full of rationals here suggest that some sort of $p$-adic analysis
might be interesting. The appearance of the square root suggests that
there is a relation to complex multiplication, or that one may have
interesting results on the field of Gaussian integers.

\section{Evaluating the Contour Integral}

The goal of this section is to relate the sum $S_{k}(x)$ to the contour
integral $I_{k}(x)$. 
\begin{thm}
For odd integers $k$, one has that \[
S_{k}(x)=I_{k}(x)+(-1)^{(k-1)/2}\left(\frac{x}{2\pi}\right)^{k-1}S_{k}\left(\frac{4\pi^{2}}{x}\right)\]
For the remainder of this text, this will be referred to as the {}``functional
equation for $S_{k}$''.\end{thm}
\begin{proof}
To prove this result, consider evaluating the contour integral $I_{s}(x)$
for a tall rectangular contour surrounding the poles at $z=1,0,-1,\ldots,1-s$.
Thus, write $I_{s}(x)=A+B+C+D$ with $A$ being the integral from
$c-ih$ to $c+ih$ for a constant $c>1$ and the height $h$ large,
eventually taking the limit $h\to\infty$. That is, $A$ forms the
right hand side of the rectangular contour. In the limit of $h\to\infty$,
one has by definition \[
A=S_{s}(x)\]

Let $B$ and $C$ be the top and bottom of the contour, so that for
$B$, the integral runs from $c+ih$ to $ih-s-\epsilon$ leftwards.
For $C$, the integral runs rightwards from $-ih-s-\epsilon$ to $c-ih$;
here we take $\epsilon>0$. The integral $D$ on the left hand side
of the rectangle closes the contour, running downwards, from $ih-s-\epsilon$
to $-ih-s-\epsilon$.

The integrals $B$ and $C$ will vanish in the limit of the height
$h\to\infty$. This can be easily seen after a simple change of variable:
\begin{eqnarray*}
B & = & \frac{1}{2\pi i}\int_{ih+c}^{ih-s-\epsilon}\frac{\Gamma(z)}{x^{z}}\zeta(s+z)\zeta(z)\, dz\\
 & = & \frac{1}{2\pi i}\int_{c}^{-s-\epsilon}\frac{\Gamma(u+ih)}{x^{u+ih}}\zeta(s+u+ih)\zeta(u+ih)\, du\end{eqnarray*}
Much of the integrand is $\mathcal{O}(1)$ in $h$, or polynomially
thereabouts. The integrand is dominated by the Gamma function, which,
from Stirling's approximation, may be seen to be \[
\Gamma(u+ih)=\mathcal{O}\left(e^{-\pi h/2}\right)\]
The complex conjugate argument applies to $C$, and thus $B$ and
$C$ vanish in the limit of $h\to\infty$. 

To evaluate $D$, begin by writing\begin{eqnarray*}
D & = & -\frac{1}{2\pi i}\int_{-ih-s-\epsilon}^{ih-s-\epsilon}\frac{\Gamma(z)}{x^{z}}\zeta(s+z)\zeta(z)\, dz\\
 & = & -\frac{1}{2\pi i}x^{s+\epsilon}\int_{-ih}^{ih}x^{u}\Gamma(-u-s-\epsilon)\zeta(-u-\epsilon)\zeta(-u-s-\epsilon)\, du\end{eqnarray*}
after a change of variable $z=-u-s-\epsilon$. One then applies the
functional equations for Gamma:\[
\Gamma(1-z)\Gamma(z)=\frac{\pi}{\sin\pi z}\]
and for zeta:\[
\zeta(s)=2^{s}\pi^{s-1}\sin\frac{\pi s}{2}\Gamma(1-s)\zeta(1-s)\]
to obtain \begin{eqnarray*}
D & =\frac{1}{2\pi^{2}i}\left(\frac{x}{2\pi}\right)^{s+\epsilon}\left(\frac{1}{2\pi}\right)^{\epsilon}\int_{-ih}^{ih}\left(\frac{x}{4\pi^{2}}\right)^{u} & \frac{\sin\frac{u+\epsilon}{2}\pi\;\sin\frac{u+s+\epsilon}{2}\pi}{\sin(u+s+\epsilon)\pi}\;\times\\
 &  & \Gamma(1+u+\epsilon)\zeta(1+u+\epsilon)\zeta(1+s+u+\epsilon)\, du\end{eqnarray*}
Another change of variable, this time as $w=1+u+\epsilon$, takes
the integral to a slightly more recognizable form: \[
D=-\frac{1}{\pi i}\left(\frac{x}{2\pi}\right)^{s-1}\int_{1+\epsilon-ih}^{1+\epsilon+ih}\left(\frac{x}{4\pi^{2}}\right)^{w}\frac{\cos\frac{w}{2}\pi\;\cos\frac{w+s}{2}\pi}{\sin(w+s)\pi}\Gamma(w)\zeta(w)\zeta(w+s)\, dw\]
Next, by taking $s=k$ to be an odd integer, the trigonometric piece
simplifies and looses its $w$ dependence: \[
\frac{\cos\frac{w}{2}\pi\;\cos\frac{w+k}{2}\pi}{\sin(w+k)\pi}=\frac{1}{2}(-1)^{(k+1)/2}\]
Pulling out this piece, one regains a recognizable integral, so that,
in the limit $h\to\infty$, one finally reaches the claimed result:\[
D=(-1)^{(k+1)/2}\left(\frac{x}{2\pi}\right)^{k-1}S_{k}\left(\frac{4\pi^{2}}{x}\right)\]
that is, \[
I_{k}(x)=A+D=S_{k}(x)-(-1)^{(k-1)/2}\left(\frac{x}{2\pi}\right)^{k-1}S_{k}\left(\frac{4\pi^{2}}{x}\right)\]
for odd integer $k$.
\end{proof}
The next section will review the application of this and the preceding
sections to specific, simple values of $x$, thus regaining many of
Plouffe's sums.

\section{Lemmas and Applications}

The first corollary demonstrates Plouffe's simplest sums for $\zeta(4m-1)$.
\begin{cor}
For $m$ integer, one has \[
I_{4m-1}(2\pi)=2S_{4m-1}(2\pi).\]
\end{cor}
\begin{proof}
Substitute $x=2\pi$ in the functional equation.
\end{proof}
This corollary provides the first concrete result of this exposition,
namely that \[
2\sum_{n=1}^{\infty}\frac{1}{n^{4m-1}\left(e^{2\pi n}-1\right)}=-\zeta(4m-1)+\frac{(2\pi)^{4m-1}}{(4m)!}\;\sum_{j=0}^{2m}\left(\begin{array}{c}
4m\\
2j\end{array}\right)\left(-1\right)^{j}B_{2j}B_{4m-2j}\]
which completely resolves one set of relationships given by Plouffe.
The functional equation opens additional possibilities. By substituting
$x=2\pi p/q$, one obtains, for $k=4m-1$, that \[
q^{k-1}\sum_{n=1}^{\infty}\frac{1}{n^{k}\left(e^{2\pi pn/q}-1\right)}+p^{k-1}\sum_{n=1}^{\infty}\frac{1}{n^{k}\left(e^{2\pi qn/p}-1\right)}=q^{k-1}I_{k}\left(\frac{2\pi p}{q}\right)\]
Thus, for example, by choosing $p=2$, $q=1$ and $k=3$, on obtains
\[
\zeta(3)=\frac{37\pi^{3}}{900}-\frac{2}{5}\sum_{n=1}^{\infty}\frac{1}{n^{3}}\left[\frac{4}{e^{\pi n}-1}+\frac{1}{e^{4\pi n}-1}\right]\]
and similarly, for $k=7$,\[
\zeta(7)=\frac{409\pi^{7}}{94500}-\frac{2}{5}\sum_{n=1}^{\infty}\frac{1}{n^{7}}\left[\frac{4}{e^{\pi n}-1}+\frac{1}{e^{4\pi n}-1}\right]\]
and one may proceed in a similar manner. There are an uncountable
infinity of such relations (since $p/q$ need not be rational). One
may take arbitrary linear combinations of these; or if one desires,
one may subtract to cancel out zeta terms, leaving behind an (uncountable)
infinity of relations for powers of $\pi$.

The next corollary shows that something more is needed for $\zeta(5)$,
$\zeta(9)$, and so on, since the most direct approach does not give
any information for such sums.
\begin{cor}
For $k$ an odd integer, one has \[
I_{k}(x)=(-1)^{(k+1)/2}\left(\frac{4\pi^{2}}{x}\right)^{k-1}I_{k}\left(\frac{4\pi^{2}}{x}\right)\]
\end{cor}
\begin{proof}
This may be proved by applying the functional equation twice in a
row. That is, it may be proved by substituting $x\to4\pi^{2}/x$ in
the functional equation and then employing the result. \end{proof}
\begin{cor}
For $m$ a positive integer, one has\[
I_{4m+1}(2\pi)=0.\]
\end{cor}
\begin{proof}
Substitute $x=2\pi$ in the preceding corollary.
\end{proof}
Following directly from the above is another corollary:
\begin{cor}
For $m$ a positive integer, one has \[
0=\sum_{j=0}^{2m+1}(-1)^{j}\frac{B_{4m-2j+2}}{(4m-2j+2)!}\frac{B_{2j}}{(2j)!}\]
\end{cor}
\begin{proof}
Write out the value of $I_{4m+1}(2\pi)$ in detail.
\end{proof}
From the above, it should be clear that the functional equation for
$S_{k}$ does not provide any statements about $S_{k}$ when $k=4m+1$.
To obtain results on sums involving $k=4m+1$, one must introduce\[
T_{s}(x)=\sum_{n=1}^{\infty}\frac{1}{n^{s}\left(e^{xn}+1\right)}\]

\begin{thm}
One has\[
T_{s}(x)=S_{s}(x)-2S_{s}(2x)\]
\end{thm}
\begin{proof}
This may be proved by re-writing in terms of the polylogarithm, along
the lines of the earlier development: \[
T_{s}(x)=-\sum_{m=1}^{\infty}(-1)^{m}{\rm Li}_{s}\left(e^{-xm}\right)\]
The even and odd terms are regrouped, as \[
T_{s}(x)=\sum_{m=1}^{\infty}{\rm Li}_{s}\left(e^{-xm}\right)-2{\rm Li}_{s}\left(e^{-2xm}\right)\]
which is seen to be a sum of $S_{s}$'s.
\end{proof}
The results for $\zeta(5)$, etc. follow from a critical observation:
that \[
S_{s}(x+2\pi i)=S_{s}(x)\]
is a periodic function. This periodicity is employed directly in the
next theorem.
\begin{thm}
For positive integer $m$, one has \begin{eqnarray*}
S_{4m+1}(2\pi) & = & I_{4m+1}(2\pi(1+i))\\
 &  & +(-1)^{m}\left[\frac{T_{4m+1}(2\pi)}{4^{m}}+2\cdot4^{m}S_{4m+1}(2\pi)-4^{m}I_{4m+1}(\pi)\right]\end{eqnarray*}
\end{thm}
\begin{proof}
Using periodicity, one writes, for $k=4m+1$, \[
S_{k}(2\pi+2\pi i)=S_{k}(2\pi)=I_{k}(2\pi(1+i))+2^{(k-1)/2}e^{3i\pi(k-1)/4}S_{k}(\pi(1-i))\]
The series $S_{k}(\pi(1-i))$ doesn't have an imaginary part; rather,
it is an alternating series, which may be expanded and written as
\[
S_{k}(\pi(1-i))=-T_{k}(\pi)+\frac{T_{k}(2\pi)+S_{k}(2\pi)}{2^{k}}\]
The $T_{k}(\pi)$ term may be eliminated by writing \[
T_{k}(\pi)=S_{k}(\pi)-2S_{k}(2\pi)\]
and the $S_{k}(\pi)$ term may be eliminated by \begin{eqnarray*}
S_{k}(\pi) & = & I_{k}(\pi)+2^{1-k}S_{k}(4\pi)\\
 & = & I_{k}(\pi)+2^{-k}\left[S_{k}(2\pi)-T_{k}(2\pi)\right]\end{eqnarray*}
Performing the various substitutions suggested above proves the theorem. 
\end{proof}
As an example of the application of the above theorem, take $m=1$,
that is, $k=5$. One easily finds that \[
I_{5}(\pi)=-\frac{15}{32}\zeta(5)+\frac{\pi^{5}}{9\cdot64}\]
and that \[
I_{5}(2\pi(1+i))=-\frac{5}{2}\zeta(5)+\frac{\pi^{5}}{9\cdot15}\]
Combining these, one gets\[
\zeta(5)=\frac{\pi^{5}}{294}-\frac{2}{35}\left[T_{5}(2\pi)+36S_{5}(2\pi)\right]\]
which is given by Plouffe. The theorem may be used to generate similar
expressions for all $\zeta(4m+1)$.

Curiously, the theorem yields results for $m=0$ as well. In this
case, one finds \[
S_{1}(2\pi)+T_{1}(2\pi)=\frac{\pi}{6}-\frac{3}{4}\log2\]

Many identities for $\pi$ are possible by taking two different expressions
for a given zeta, and subtracting them, leaving behind a rational
combination of the sums and $\pi$. Thus, for example, the following
theorem for Apéry's constant:
\begin{thm}
A series expression for $\pi^{3}$ is given by \[
\pi^{3}=720\cdot S_{3}(\pi)-900\cdot S_{3}(2\pi)+180\cdot S_{3}(4\pi)\]
 \end{thm}
\begin{proof}
This follows by taking the general expression for $k=4m-1$:\begin{eqnarray*}
16^{m}S_{4m-1}(\pi)+4S_{4m-1}(4\pi) & = & 16^{m}I_{4m-1}(\pi)\\
 & = & -\frac{1}{2}\zeta(4m-1)\left[16^{m}+4\right]\\
 &  & -(2\pi)^{4m-1}\sum_{j=0}^{2m}(-4)^{j}\frac{B_{4m-2j}}{(4m-2j)!}\frac{B_{2j}}{(2j)!}\end{eqnarray*}
solving for $\zeta(4m-1)$ and then using \[
2S_{4m-1}(2\pi)=I_{4m-1}(2\pi)=-\zeta(4m-1)+(2\pi)^{4m-1}\;\sum_{j=0}^{2m}\left(-1\right)^{j}\frac{B_{4m-2j}}{(4m-2j)!}\frac{B_{2j}}{(2j)!}\]
to eliminate the appearance of the zeta. The resulting expression
may then be solved for $\pi^{4m-1}$.
\end{proof}
Relationships involving square roots also arise naturally. 
\begin{thm}
One has \[
\zeta(3)=\frac{5\pi^{3}}{72}-\frac{1}{2}\, S_{3}\left(2\pi\sqrt{3}\right)-\frac{3}{2}\, S_{3}\left(\frac{2\pi\sqrt{3}}{3}\right)\]
\end{thm}
\begin{proof}
This follows from the general expression given earlier, that \[
S_{k}\left(\frac{2\pi p}{q}\right)+\left(\frac{p}{q}\right)^{k-1}S_{k}\left(\frac{2\pi q}{p}\right)=I_{k}\left(\frac{2\pi p}{q}\right)\]
Here, making the substitution $q=\sqrt{p}$ one obtains \begin{eqnarray*}
S_{k}\left(2\pi\sqrt{p}\right)+p^{(k-1)/2}S_{k}\left(\frac{2\pi\sqrt{p}}{p}\right) & = & -\frac{1}{2}\zeta(k)\left[1-(-p)^{(k-1)/2}\right]\\
 &  & +\frac{(-1)^{(k-1)/2}}{2\sqrt{p}}(2\pi)^{k}\;\\
 &  & \;\;\sum_{j=0}^{(k+1)/2}\left(-p\right)^{j}\frac{B_{k+1-2j}}{(k+1-2j)!}\frac{B_{2j}}{(2j)!}\end{eqnarray*}
The specific result follows after choosing $k=3$ and $p=3$.
\end{proof}

\section{Some Bernoulli number identities}

In addition to the previously noted identity \[
0=\sum_{j=0}^{2m+1}(-1)^{j}\frac{B_{4m-2j+2}}{(4m-2j+2)!}\frac{B_{2j}}{(2j)!}\]
there are several other identities on sums of Bernoulli numbers that
result from the previous developments. These are briefly stated here. 
\begin{thm}
For integer $m$, one has\[
0=\sum_{j=0}^{m}(-4)^{j}\left[\frac{B_{4m-4j+2}}{(4m-4j+2)!}\frac{B_{4j}}{(4j)!}+2\frac{B_{4m-4j}}{(4m-4j)!}\frac{B_{4j+2}}{(4j+2)!}\right]\]
\end{thm}
\begin{proof}
Consider the sums resulting from $0=I_{4m+1}(2\pi(1+i))-I_{4m+1}(2\pi(1-i))$.\end{proof}
\begin{thm}
For integer $m$, one has \begin{eqnarray*}
\sum_{k=0}^{2m}(-1)^{k}\frac{B_{4m-2k}}{(4m-2k)!}\frac{B_{2k}}{(2k)!} & = & \sum_{j=0}^{m}(-4)^{j}\frac{B_{4m-4j}}{(4m-4j)!}\frac{B_{4j}}{(4j)!}\\
 &  & -2\sum_{j=0}^{m-1}(-4)^{j}\frac{B_{4m-4j-2}}{(4m-4j-2)!}\frac{B_{4j+2}}{(4j+2)!}\end{eqnarray*}
 \end{thm}
\begin{proof}
Consider the sums resulting from the identity \[
I_{4m-1}(2\pi)=I_{4m-1}(2\pi(1+i))+I_{4m-1}(2\pi(1-i))\]

\end{proof}

\section{Modular Relations}

The various sums and quantities above can be seen to be quasi-modular
by making a simple change of variable, namely by making the substitution
$x=2\pi i\tau$. By {}``quasi-modular'', it is meant that the various
terms almost have simple behaviors under the Mobius transformation
$\tau\to(a\tau+b)/(c\tau+d)$ for integer $a,b,c$ and $d$. In this
respect, the sums bear close resemblance to theta functions, which
obey similar relations. These relationships are brought to focus here. 

First, define $K_{k}(\tau)=I_{k}(2\pi i\tau)$. This change of variable
results in a definition which seems simpler than that for $I_{k}$.
An expansion in $\tau$ is often referred to as a {}``Fourier series''
in the context of hyperbolic geometry: \[
K_{k}(\tau)=\frac{\tau^{k-1}-1}{2}\zeta(k)-\frac{(2\pi i)^{k}}{2\tau}\;\sum_{j=0}^{(k+1)/2}\tau^{2j}\frac{B_{2j}}{(2j)!}\frac{B_{k+1-2j}}{(k+1-2j)!}\]
This quantity is almost a modular form of weight $k-1$, in that \[
K_{k}\left(\frac{-1}{\tau}\right)=-\tau^{1-k}K_{k}(\tau)\]
Its only {}``almost'' a modular form, because it is not periodic
in $\tau$, that is \[
K_{k}(\tau+1)\ne K_{k}(\tau)\]
By contrast, $P_{k}(\tau)=S_{k}(2\pi i\tau)$ is periodic: \[
P_{k}(\tau+1)=P_{k}(\tau)\]
but is not quite modular under inversion: \[
P_{k}\left(\frac{-1}{\tau}\right)=\tau^{1-k}P_{k}(\tau)+K_{k}\left(\frac{-1}{\tau}\right)\]
One may define a simple variant of the sums that does have a simple
transformation under inversion, namely \[
M_{k}(\tau)=P_{k}(\tau)-\frac{1}{2}K_{k}(\tau)\]
which transforms as \[
M_{k}\left(\frac{-1}{\tau}\right)=-\tau^{1-k}M_{k}(\tau)\]
However, $M_{k}$ is then not periodic. 

\begin{figure}
\caption{Phase of $P_{-1}(q)$ on the unit disk}

\includegraphics[scale=0.6,bb=0in 0in 601pt 601pt]{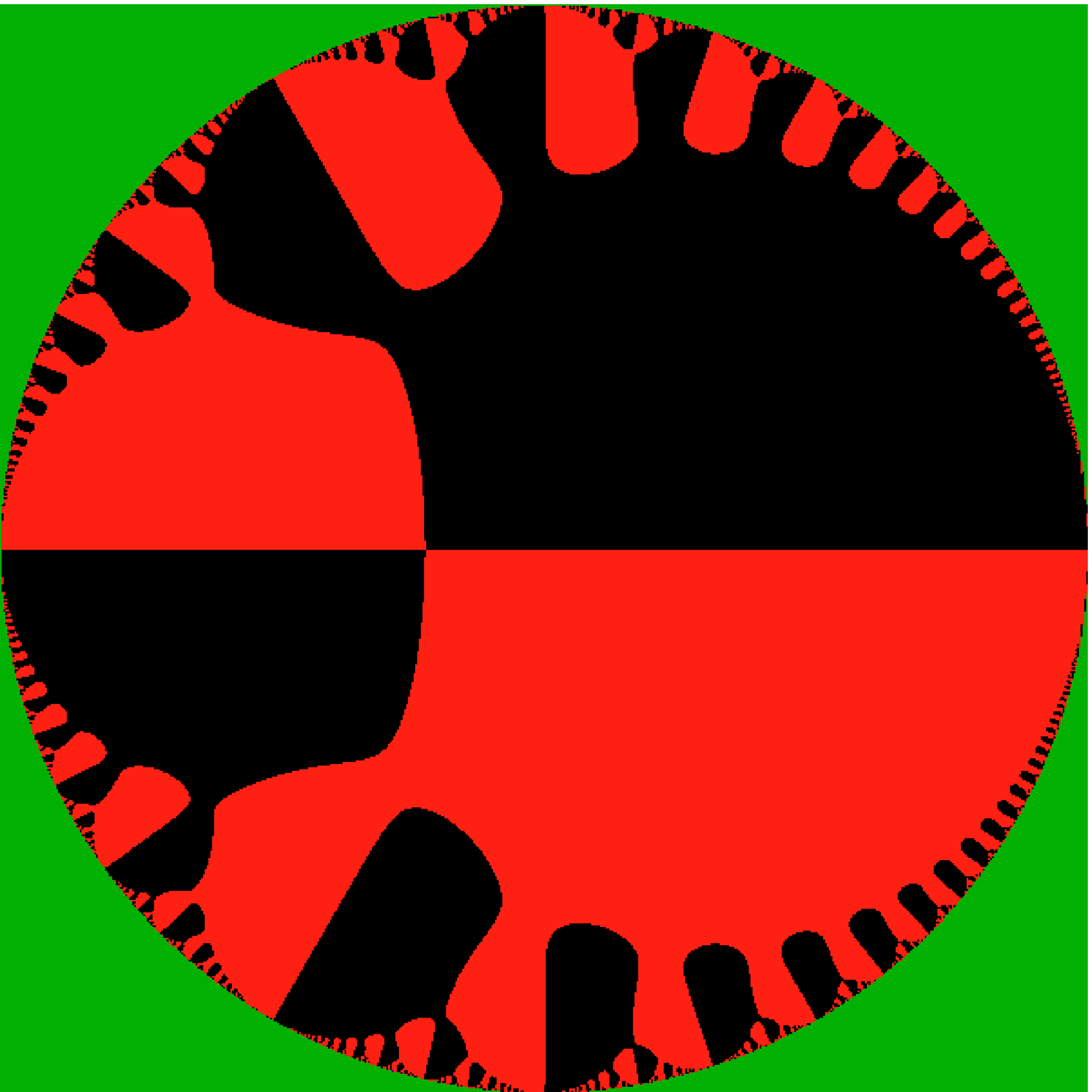}

This graphic shows the phase\[
\arg P_{-1}(q)=\arg\sum_{n=1}^{\infty}\frac{n}{q^{n}-1}\]
 where $\arg f(z)=\Im\log f(z)$ is the usual arg of a function. The
color scheme is such that black represents areas where $\arg>0$ and
red represents areas where $\arg<$0. The absence of other colors
indicates that the phase is rather closely confined to the vicinity
of $0$ for most all of the disk. Numerically, the absolute value
of the phase is smaller than $10^{-3}$ for much of the disk. In particular,
this indicates that there are no zeros at all in the interior of the
disk, as a zero would be surrounded by a region where the phase wraps
around by $2\pi$. 

The function $P_{s}(q)$ does have poles at $q=e^{2\pi im/n}$ for
all rationals $m/n$; these are visible at the edges of the disk.
The fractal nature of this image is the characteristic signature of
a modular form of weight 2; the self-similar regions are just copies
of the fundamental region of the modular group $SL(2,\mathbb{Z})$.

\lyxline{\normalsize}
\end{figure}

\begin{figure}
\caption{Phase of $P_{-5}(q)$ on the unit disk}

\includegraphics[clip,scale=0.6,bb=0in 0in 601pt 601pt]{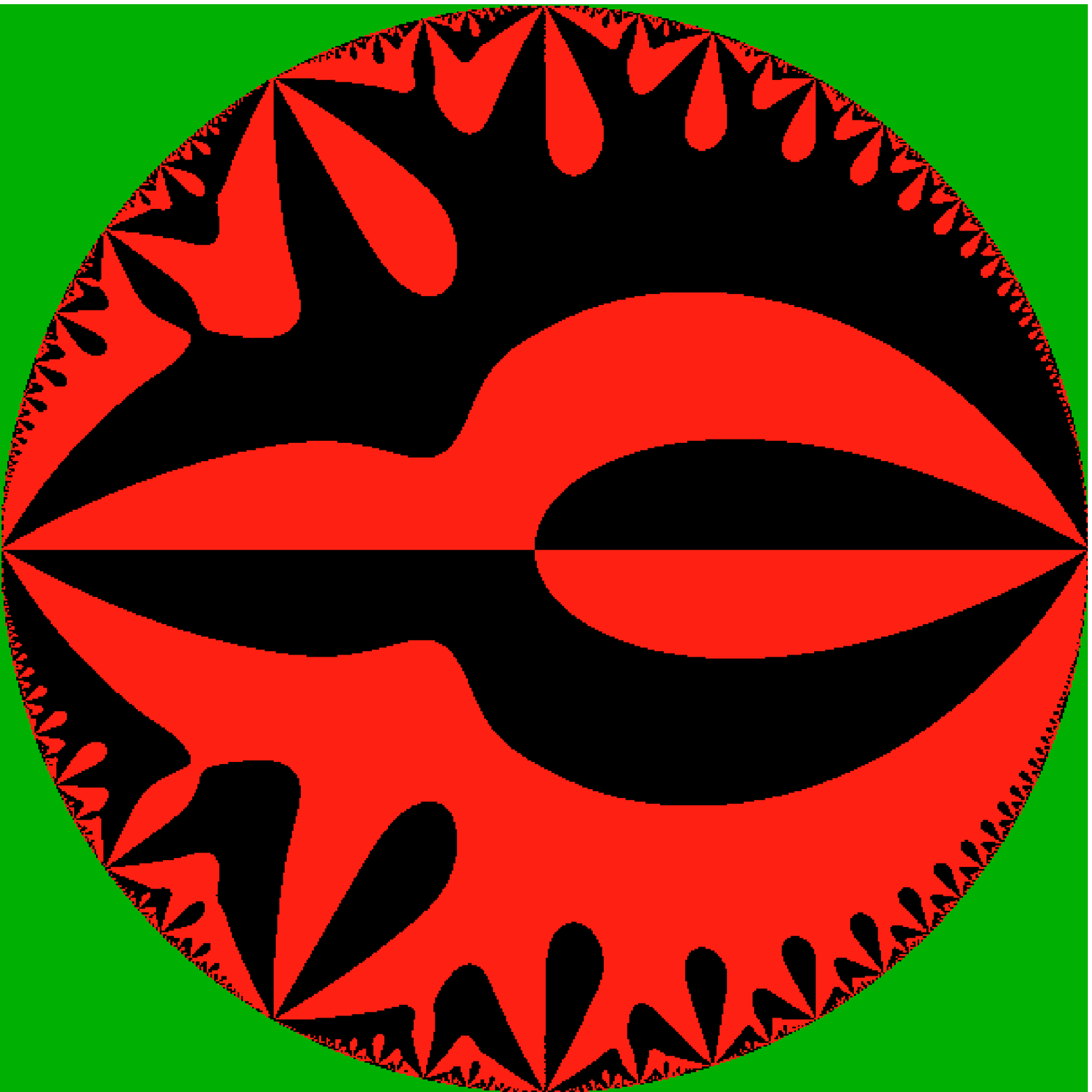}

This graphic shows the phase\[
\arg P_{-5}(q)=\arg\sum_{n=1}^{\infty}\frac{n^{5}}{q^{n}-1}\]
The color scheme is as in the previous image, and similar conclusions
apply: there are no zeros at all in the interior of the disk. The
absolute value of the phase is tiny: numerically, it is within $10^{-9}$
of zero for much of the disk.

The fractal nature of this image is the characteristic signature of
a modular form of weight 6; it should be compared to the image of
the modular invariant $g_{3}$ shown in the next image. 

\lyxline{\normalsize}
\end{figure}

\begin{figure}
\caption{Graph of the modular invariant $g_{3}$.}

\includegraphics[scale=0.6,bb=0in 0in 701pt 701pt]{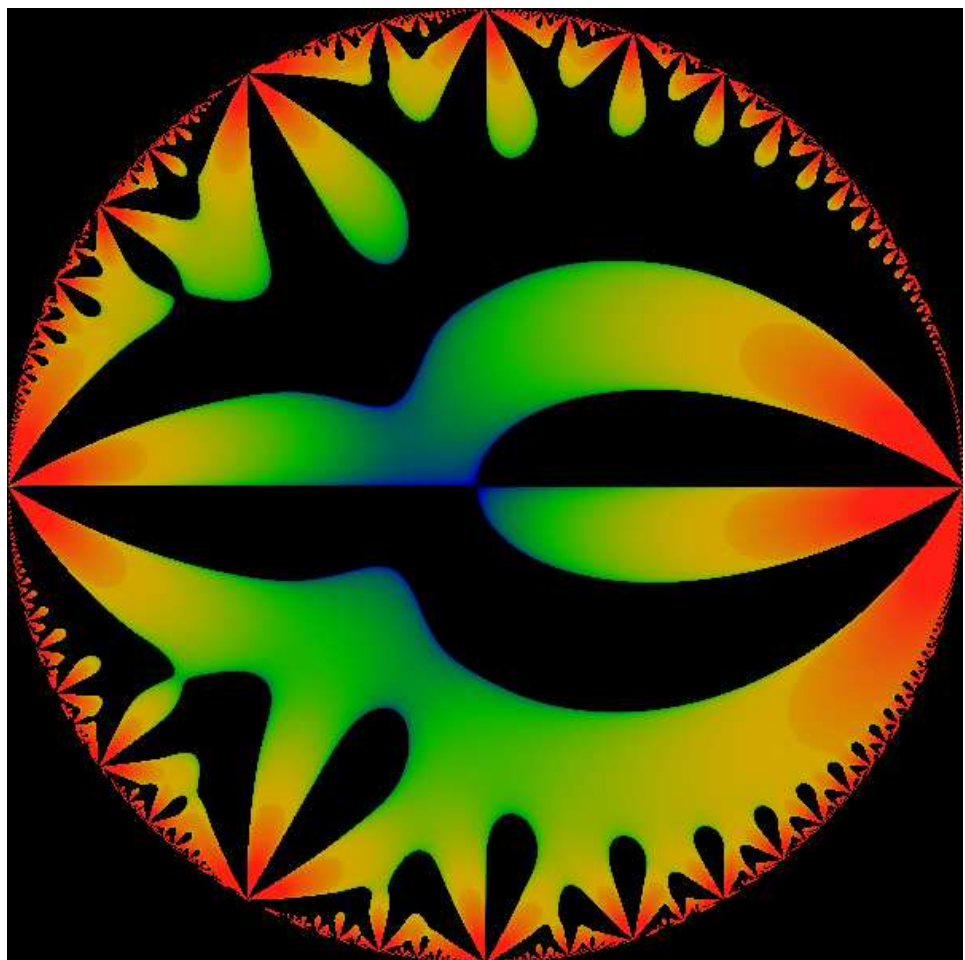}

This figure shows the imaginary part of the modular invariant $g_{3}$,
one of the invariants of an elliptic curve. Specifically, it shows
the imaginary part of \[
g_{3}(q)=\frac{8\pi^{6}}{27}\left[1-504\sum_{n=1}^{\infty}\frac{n^{5}q^{n}}{1-q^{n}}\right]\]
 which is a modular form of weight 6. The colors are chosen such that
black represents areas that are negative, blue and green represent
areas with smaller values, and red those areas with the largest values.

\lyxline{\normalsize}
\end{figure}

The analytic structure of $P_{k}$ is curious: it has a pole at $\tau=i\infty$
and. by the inversion formula and periodicity, at every rational value
of $\tau$. This is clearly visible in the graphic below, which shows
$P_{3}$ on the $q$-series or {}``punctured disk'' coordinates
$q=e^{2\pi i\tau}$:\[
P_{s}(q)=\sum_{n=1}^{\infty}\frac{1}{n^{s}\left(q^{n}-1\right)}\]

In this form, the relation to modular forms becomes the most immediate.
Consider the Lambert series:\[
\sum_{n=1}^{\infty}\frac{1}{n^{s}}\frac{q^{n}}{1-q^{n}}=-\zeta(s)-P_{s}(q)=\sum_{m=1}^{\infty}\sigma_{-s}(m)q^{m}\]
Here, $\sigma_{s}(m)$ is the divisor function: \[
\sigma_{s}(m)=\sum_{n|m}\; n^{s}\]
with the notation $n|m$ denoting that the sum extends over all divisors
$n$ of $m$. This should be compared to the Eisenstein series \cite[see section 3.10]{Apo90}
\[
G_{2k}(q)=2\zeta(2k)+\frac{2(2\pi i)^{2k}}{(2k-1)!}\sum_{m=1}^{\infty}\sigma_{2k-1}(m)q^{m}\]
That is, the sums $P_{s}$ can be re-written in terms of the Eisenstein
series $G_{1-s}$, and the behaviour of one under the action of the
modular group can be igven in terms of the other.

\section{Conclusions}

Most of the sums discussed here are suggestive of linear algebra.
So for example, one may write the sum $S_{k}$ as the dot-product
between the (infinite-dimensional) vector $n^{-k}$ and the vector
$\left(e^{2\pi n}-1\right)^{-1}$. The evaluation of these sums forms
a regular pattern in $k$, suggesting that, for example, that $n^{-k}$
could be taken to be the matrix elements of some linear operator.
However, the significance of this operator (aside from assorted shallow
results and relations) is completely unclear. The sums over the Bernoulli
numbers are even reminiscent of some crazy Atiyah-Singer-like indexing;
but the underlying operators are utterly unclear. Put another way,
it is well-known in physics and mathematics that regular patterns
are the result of symmetries; the sums discussed here form a regular
pattern, but the nature of the symmetry that generates it is unclear.
During manipulations, one gets the sense that there are plenty of
other relations to be discovered; certainly, the mere heft of Ramanujan's
tomes suggest as much. Yet, there is no picture of a generator that
can be operated to generate the myriad of relations; some symmetry
group presentation seems to be missing. A similar problem exists in
the theory of hypergeometric series, where there is an embarrassment
of riches in terms of relations and identities, and yet a unifying
theory is lacking. The ingredients to the sums discussed here include
the Gamma function and the Bernoulli polynomials; among many other
properties, these have a common set of relations in the $p$-adic
{}``multiplication theorems''; the sums here vaguely resemble the
multiplication theorems of characteristic zero. There are also similar
phenomena and sums that occur in the theory of dynamical systems,
and in particular, in symbolic dynamics; there, a group structure,
or at least, a monoid structure, together with an explicit treatment
in terms of linear operators, is more common. Any of these connections
present intriguing avenues for future research; however, the overall
problem, of discovering the underlying symmetry that leads to such
relations, seems unattainably hard to solve.

Thanks to Simon Plouffe for generating interest in such sums.

\section{Appendix: Polylog Integral}

This appendix proves the following theorem:
\begin{thm}
The polylogarithm may be written as the integral \[
{\rm Li}_{s}\left(e^{-u}\right)=\frac{1}{2\pi i}\int_{c-i\infty}^{c+i\infty}\Gamma(r)\zeta(r+s)u^{-r}dr\]
\end{thm}
\begin{proof}
The proof below is cribbed from the Wikipedia article on polylogarithms\cite{WP-PL}.
One begins by writing the the Mellin transform of the polylog, as
\[
M_{s}(r)=\int_{0}^{\infty}{\rm Li}_{s}\left(ye^{-u}\right)u^{r}\frac{du}{u}\]
Using an integral representation of the polylogarithm, \[
{\rm Li}_{s}\left(w\right)=\frac{1}{\Gamma(s)}\int_{0}^{\infty}\frac{t^{s-1}}{w^{-1}e^{t}-1}dt\]
 and substituting, one obtains \[
M_{s}(r)=\frac{1}{\Gamma(s)}\int_{0}^{\infty}\int_{0}^{\infty}\frac{u^{r-1}t^{s-1}}{y^{-1}e^{t+u}-1}dt\, du\]
A change of variable $t=ab$ and $u=a(1-b)$ with $dt\, du=a\, da\, db$
gives \begin{eqnarray*}
M_{s}(r) & = & \frac{1}{\Gamma(s)}\int_{0}^{1}b^{s-1}(1-b)^{r-1}db\int_{0}^{\infty}\frac{a^{r+s-1}}{y^{-1}e^{a}-1}\, da\\
 & = & \Gamma(r)\,{\rm Li}_{r+s}(y)\end{eqnarray*}
The inverse Mellin transform may now be employed to write \[
{\rm Li}_{s}\left(ye^{-u}\right)=\frac{1}{2\pi i}\int_{c-i\infty}^{c+i\infty}u^{-r}\Gamma(r)\,{\rm Li}_{r+s}(y)\, dr\]
By setting $y=1$, one then uses ${\rm Li}_{s+r}(1)=\zeta(s+r)$ to
obtain the desired result.
\end{proof}
\bibliographystyle{plain}
\bibliography{/home/linas/linas/fractal/paper/fractal}

--

--
\end{document}